\documentclass[11pt,a4paper]{amsart}
\usepackage[T1]{fontenc}
\usepackage[utf8]{inputenc}
\usepackage[english]{babel}
\usepackage{lmodern}
\usepackage{microtype}
\usepackage{csquotes}
\usepackage[a4paper,textwidth=145mm,textheight=230mm,centering]{geometry}
\usepackage{amssymb,amsmath,amsthm}
\usepackage[backend=biber,style=numeric,sorting=nyt,giveninits=true,isbn=false,url=false,doi=true,eprint=true]{biblatex}
\usepackage[hidelinks]{hyperref}
\hypersetup{pdftitle={Failure of endpoint integrability for orthogonal projections of Frostman measures},pdfsubject={Geometric measure theory},pdfkeywords={Orthogonal projections, Frostman measures, endpoint integrability, Lorentz spaces, Hausdorff measure}}
\newtheorem{theorem}{Theorem}
\newtheorem{lemma}[theorem]{Lemma}
\title[Failure of endpoint integrability for projections]{Failure of endpoint integrability for orthogonal projections of Frostman measures}
\author{Vicente Vergara}
\address{Department of Mathematics, Faculty of Physical and Mathematical Sciences, University of Concepci\'on, Concepci\'on, Chile.}
\email{vvergaraa@udec.cl}
\date{}

\begin{document}
\raggedbottom
\begin{abstract}
For compactly supported $s$-Frostman measures in the plane, with $1<s<2$, almost all orthogonal projections have densities in $L^p$ for $1\le p<2/(2-s)$. Orponen established the sharpness of this threshold, leaving open whether every such measure has even a single projection with a density in the endpoint space. We answer this question in the negative, even for weak endpoint integrability. For every $1<s<2$, we construct a compact set $E\subset[0,1]^2$ with $0<\mathcal{H}^s(E)<\infty$ such that $\mathcal{H}^s|_E$ is $s$-Frostman and none of its orthogonal projections admits a density in $L^{2/(2-s),\infty}$. More precisely, in every direction the projected mass of an interval $J$, normalised by $|J|^{s/2}$, is unbounded as $J$ varies. The construction distributes the orientations among different parent squares and generations while preserving the exact Frostman exponent.
\end{abstract}
\maketitle

\section{Introduction}
Marstrand's theorem \cite{Marstrand1954} states that if a planar Borel set has Hausdorff dimension greater than $1$, almost all of its orthogonal projections onto lines have positive Lebesgue measure. This conclusion does not, however, describe how mass on the original set is distributed under projection. Kaufman's energy approach \cite{Kaufman1968} and its development through Fourier analysis (see \cite{Mattila2015}) provide a way to study the regularity of projected measures from a dimensional bound on mass.

Let $\mu$ be a compactly supported finite positive Borel measure on $\mathbb{R}^2$. We say that $\mu$ is $s$-Frostman if there exists $C_\mu>0$ such that
\begin{equation}\label{eq:frostman}
 \mu(B(x,r))\le C_\mu r^s,\qquad x\in\mathbb{R}^2,\quad r>0.
\end{equation}
For $e\in S^1$, write $\pi_e(x)=x\cdot e$ and $\nu_e=(\pi_e)_\#\mu$. Let $\sigma$ denote normalised arc length measure on $S^1$. We write $\nu_e\in L^p$ if $\nu_e$ is absolutely continuous and its density $f_e$ belongs to $L^p(\mathbb{R})$, and adopt the same convention for the weak space $L^{p,\infty}$.

Fix $1<s<2$. Condition \eqref{eq:frostman} ensures that the energy $I_t(\mu)=\iint |x-y|^{-t}\,\mathrm{d}\mu(x)\,\mathrm{d}\mu(y)$ is finite for $0<t<s$. By the Fourier slice identity and the polar-coordinate representation of the energy, for $1<t<s$ the densities $f_e$ belong to $H^{(t-1)/2}(\mathbb{R})$ for almost every $e$, and their squared Sobolev norms are integrable with respect to $\sigma$. The one-dimensional Sobolev embedding and compact support give (see \cite[Section~3.1]{DabrowskiOrponenVilla2022})
\begin{equation}\label{eq:subcritical}
 \int_{S^1}\|f_e\|_{L^p(\mathbb{R})}^2\,\mathrm{d}\sigma(e)<\infty,
 \qquad 1\le p<p_*(s):=\frac{2}{2-s}.
\end{equation}
Finiteness of the critical energy $I_s(\mu)$ would give integrability at $p_*(s)$ by the same argument, but is not guaranteed by \eqref{eq:frostman}.

Peres and Schlag \cite[Section~9.2(ii)]{PeresSchlag2000} asked whether the threshold $p_*(s)$ could be improved. Orponen \cite[Theorem~1.1]{Orponen2025} answered this question in the negative: for every $p>p_*(s)$ he constructed an $s$-Frostman measure none of whose projections belongs to $L^p$. Whether the endpoint itself is admissible, even in a single direction for each measure, remained open (\cite[Section~1.1]{Orponen2025}).

This question has a geometric interpretation. The mass $\nu_e(J)$ of an interval $J$ is the mass of the strip $\pi_e^{-1}(J)$ of width $|J|$. H\"older's inequality gives, for $f_e\in L^p$ with $1<p<\infty$,
\[
 \nu_e(J)\le \|f_e\|_{L^p}|J|^{1-1/p}.
\]
At the endpoint $p=p_*(s)$, the width appears with exponent $s/2$; a bound with this same power is necessary even under the weaker assumption $f_e\in L^{p_*(s),\infty}$, by \eqref{eq:weak}. We show that a single Frostman measure can violate this bound in every direction, even when it is the restriction of Hausdorff measure to an $s$-dimensional set.

\begin{theorem}\label{thm:main}
For every $1<s<2$, there exists a compact set $E\subset[0,1]^2$ with $0<\mathcal{H}^s(E)<\infty$ such that $\lambda=\mathcal{H}^s|_E$ is $s$-Frostman and
\begin{equation}\label{eq:obstruction}
 \sup_J\frac{(\pi_e)_\#\lambda(J)}{|J|^{s/2}}=\infty
 \qquad\text{for every }e\in S^1.
\end{equation}
The supremum is taken over bounded intervals of positive length. In particular,
\[
  (\pi_e)_\#\lambda\notin L^{2/(2-s),\infty}(\mathbb{R})
  \qquad\text{for every }e\in S^1.
\]
The same conclusions hold for the probability measure $\lambda/\mathcal{H}^s(E)$.
\end{theorem}

The theorem and \eqref{eq:subcritical} give the exact universal range of spatial integrability: $1\le p<p_*(s)$, whether one requires the property in almost every direction or merely in some direction. Applying \eqref{eq:subcritical} to a sequence of exponents increasing to $p_*(s)$, we see that almost all projections of the example have densities in all the subcritical spaces. Thus the endpoint obstruction does not require singular projections and is realised by Hausdorff measure on the support itself. We do not claim that the set is Ahlfors regular.

Control of the norms as the direction varies is a separate question. Dąbrowski, Orponen and Villa \cite[Thm.~1.2, Ex.~1.3 and Sec.~1.1]{DabrowskiOrponenVilla2022} proved that $\int_{S^1}\|f_e\|_{L^p}^p\,\mathrm{d}\sigma(e)<\infty$ for $1\le p<(3-s)/(2-s)$, and showed that the supremum of this range cannot be increased. They also applied these joint estimates to Furstenberg sets. Liu \cite[Proposition~10.1]{Liu2024} obtained bounds with distinct spatial and angular exponents $p$ and $q$. For $2\le q<p$, reaching $p=p_*(s)$ in these bounds requires energies of order at least $s$ when the potential orders are those guaranteed by the Frostman condition. Li and Liu \cite[Theorem~1.1]{LiLiu2026} extended the diagonal estimates to Frostman measures on the space of directions and obtained improvements under additional structural hypotheses. Our theorem settles spatial integrability at the endpoint under the Frostman condition alone, rather than the full range of mixed norms.

\subsection*{Idea of the proof}
The difficulty of constructing an endpoint example can be seen from Orponen's blocks \cite[Proposition~2.3 and Sections~2.2--2.5]{Orponen2025}. At scale $\delta$, his blocks consist of approximately $\delta^{-s}$ squares arranged along segments of length comparable to $\delta^{1-s/2}$. For each direction, one of the segments contributes a mass of order $\delta^{s/2}$ to an interval of length comparable to $\delta$. The resulting lower bound for the $L^p$ norm, by H\"older's inequality, is of order $\delta^{s/2-1+1/p}$: it grows without bound if $p>p_*(s)$, but its exponent vanishes at the endpoint.

When a block is repeated inside a parent square of side length $R$ and mass $m_Q$, the critical ratio corresponding to the mass contributed by a single segment becomes of order
\[
 \frac{m_Q\delta^{s/2}}{(R\delta)^{s/2}}=\frac{m_Q}{R^{s/2}}\le C_s C_\mu R^{s/2},
\]
where $C_s$ depends only on $s$. Decreasing the relative scale $\delta$ does not compensate for this loss. This loss pertains only to the contribution of a single segment and does not imply an upper bound on all possible concentrations arising from the block. The supercritical argument in \cite[Section~3]{Orponen2025} relies on a strict separation between auxiliary dimension exponents, which is unavailable at the endpoint.

Rather than requiring each parent square to reproduce all orientations, we distribute the angular coverage among the parents. If there are $K$ parent squares of side length $R$ and mass $1/K$, with $KR^s\asymp1$, we choose $m$ segments of $k$ squares in each parent, where, for a target amplitude $A\ge1$,
\[
 k\asymp A\sqrt K\,\delta^{-s/2},\qquad mk\asymp\delta^{-s}.
\]
The implicit constants in these comparisons may depend on $s$ and on the uniform bounds for $KR^s$. Each segment covers an angular interval of length of order $1/k$. The $Km$ segments, with orientations interleaved among the parents, jointly cover a sector of length of order $Km/k\asymp A^{-2}$. For each projection whose kernel direction lies in the sector, the contribution of one segment has a critical ratio of order
\[
 \frac{1}{Km(R\delta)^{s/2}}\asymp\frac{A}{\sqrt{KR^s}}\asymp A.
\]
The resulting ratio grows with $A$, with no loss of a power of $\delta$.

Two estimates make it possible to iterate this idea. The lemma in the next section ensures uniform packing even in very narrow sectors by taking $\delta$ sufficiently small after fixing the angular width. Adjusting the cardinalities with summable multiplicative errors then preserves $1\le K_nR_n^s\le\exp(1)$ at every generation; this invariant controls the Frostman condition and the Hausdorff measure of the limit set. For each value in a sequence $A\to\infty$, we cover the space of directions by finitely many sectors, one per generation. The nesting ensures that the mass in each of these concentrations is preserved. In this way, we obtain unbounded critical ratios in all directions without requiring the corresponding concentrations to occur simultaneously at the same scale.

Throughout the proof, $s$ is fixed. We use the notation $X\lesssim Y$ and $X=O(Y)$ to mean, respectively, $X\le CY$ and $|X|\le CY$, with a constant that may depend on $s$ but not on the scales, cardinalities or angular parameters. We write $X\asymp Y$ when both corresponding inequalities hold; a subscript indicates the explicit dependencies of the constant. Generic constants may vary from line to line. Dyadic squares are closed and are chosen to be pairwise disjoint. Segment directions are considered modulo $\pi$; angles in a sector are represented by real numbers in an interval of length at most $\pi$.

\section{A packing lemma for directions in a sector}
\begin{lemma}\label{lem:packing}
Fix $1<s<2$. There are constants $c_0,c_1,C>0$, which may depend on $s$, with the following property. Let $0<w\le\pi$, $\delta\in 2^{-\mathbb{N}}$ and let $k,m$ be positive integers. Suppose that $L=4k\delta$ and that
\begin{equation}\label{eq:packing-cond}
  1\le mk\delta^s\le2,\quad L\le c_0,\quad
  L^{2-s}\le w,\quad \delta^{2-s}\le c_1w.
\end{equation}
Given the angles $\theta_i=\theta_0+iw/m$, $0\le i<m$, there exists a family $\mathcal{P}$ of exactly $mk$ pairwise disjoint closed dyadic squares of side length $\delta$, contained in the interior of $[0,1]^2$ and partitioned into $m$ groups of $k$ squares. The squares in group $i$ lie in the $\sqrt2\delta$-neighbourhood of a segment of length $L$ and direction $\theta_i$. Moreover,
\begin{equation}\label{eq:packing}
 \#\{Q\in\mathcal{P}:Q\cap B(x,r)\ne\varnothing\}
 \le C(r/\delta)^s,\qquad \delta\le r\le1,
\end{equation}
for every $x\in\mathbb{R}^2$.
\end{lemma}

\begin{proof}
Fix a sufficiently small absolute constant $b>0$, for instance $b\le 1/6400$, and set
\[
 h=\min\{m,\lfloor b/L^2\rfloor\}.
\]
Choose $c_0$ sufficiently small that $b/L^2\ge2$ and $3L+\sqrt2\delta<1/4$. Using a grid, choose $h$ centres in $[1/4,3/4]^2$ with mutual separation at least $1/(4\sqrt h)$. Since $h\le b/L^2$, the separation is at least $20L$. Assign the indices $i$ to the centres according to their residue classes modulo $h$, so that each centre receives at most $\lceil m/h\rceil$ indices.

For an index $i$ assigned to a centre $c$, let $u_i=(\cos\theta_i,\sin\theta_i)$. Choose the $k$ points
\begin{equation}\label{eq:points}
  c+(2L+4\ell\delta)u_i,\qquad \ell=0,\ldots,k-1.
\end{equation}
For each point, select the dyadic square of side length $\delta$ that contains it, using half-open dyadic squares to assign points on the boundary uniquely and then taking the closure. All the squares lie in the interior of the unit square. Those assigned to a given index lie in the $\sqrt2\delta$-neighbourhood of the segment joining $c+2Lu_i$ and $c+3Lu_i$.

We first verify that the squares are disjoint. On a given ray, the points in \eqref{eq:points} are at least $4\delta$ apart. Points associated with distinct centres are at least $14L$ apart. If $h=m$, there is a single ray per centre. If $h<m$, the angles of two distinct rays at the same centre differ by at least $wh/m$ and at most $\pi$. For unit vectors $u,v$ forming an angle $d\in[0,\pi]$ and radii $\rho,\rho'\ge2L$, the identity
\[
 |\rho u-\rho'v|^2=(\rho-\rho')^2+2\rho\rho'(1-\cos d)
\]
implies $|\rho u-\rho'v|\gtrsim Ld$, since $2\sin(d/2)\ge2d/\pi$. As $h\ge b/(2L^2)$ and $mL=4mk\delta\le8\delta^{1-s}$, the distance between points on distinct rays at the same centre is bounded below by
\[
  c'\frac{Lwh}{m}\gtrsim \frac{w}{mL}\gtrsim w\delta^{s-1},
\]
where $c'>0$ is absolute and the constants may depend on the fixed absolute constant $b$. By the last condition in \eqref{eq:packing-cond}, taking $c_1$ sufficiently small makes this distance greater than $2\sqrt2\delta$. The preceding separation bounds ensure that no two selected closed squares intersect, since every point of each square lies at distance at most $\sqrt2\delta$ from its sampling point. Hence there are exactly $mk$ squares.

We now prove \eqref{eq:packing}. For $\delta\le r\le L$, a ball of radius $r$ can meet squares associated with only one centre, by the separation of the centres and the upper bound $3L+\sqrt2\delta$ for the distance from a square to its centre. If the ball meets squares on several rays, the corresponding sampling points lie in a ball of radius $r+\sqrt2\delta$. The same distance identity implies that the corresponding angles lie in an interval of length $O(r/L)$. Each ray contributes $O(r/\delta)$ squares. When $h=m$, there is only one ray per centre and the required bound follows. When $h<m$, the angular spacing gives at most $O(1+rm/(Lwh))$ rays. Since $m/(Lh)\lesssim mL\lesssim\delta^{1-s}$, we obtain
\begin{align*}
 \#\{Q\in\mathcal{P}:Q\cap B(x,r)\ne\varnothing\}
 &\lesssim \frac r\delta+\frac{r^2}{w\delta^s}\\
 &\lesssim (r/\delta)^s.
\end{align*}
For the last inequality we use $r\ge\delta$, $s>1$ and
\[
  \frac{r^2}{w\delta^s}
  =(r/\delta)^s\frac{r^{2-s}}w\le(r/\delta)^s,
\]
since $r^{2-s}\le L^{2-s}\le w$.

For $L\le r\le1$, the separation of the centres and the fact that the squares lie at distance $O(L)$ from their centres imply that a ball meets squares associated with $O(1+hr^2)$ centres. Each centre has at most $\lceil m/h\rceil k\le2mk/h$ squares. The definition of $h$ gives $1/h\lesssim 1/m+L^2$. Hence
\begin{align*}
 \#\{Q\in\mathcal{P}:Q\cap B(x,r)\ne\varnothing\}
 &\lesssim mk(1/h+r^2)\\
 &\lesssim k+mk(L^2+r^2)\\
 &\lesssim r/\delta+r^2\delta^{-s}\\
 &\lesssim (r/\delta)^s.
\end{align*}
Here we use $k=L/(4\delta)$, $mk\le2\delta^{-s}$ and $L\le r\le1$. This proves the lemma with constants independent of $w,m,k,\delta$.
\end{proof}

\section{A refinement step}
Suppose that a generation consists of $K$ disjoint dyadic squares of side length $R$ contained in the unit square, each assigned mass $1/K$, and that, for some constant $C_0\ge1$,
\begin{equation}\label{eq:invariant0}
  1\le KR^s\le C_0.
\end{equation}
Fix an amplitude $A\ge1$, a sector $[\theta_0,\theta_0+w]$ with $0<w\le\pi/A^2$ and an error tolerance $0<\varepsilon\le1$. We refine the $K$ squares to produce a large concentration of projected mass for every projection whose kernel direction lies in the sector.

Choose a dyadic scale $\delta>0$ and set
\begin{equation}\label{eq:km}
  k=\lceil A\sqrt K\,\delta^{-s/2}\rceil,
  \qquad m=\lceil\delta^{-s}/k\rceil.
\end{equation}
For $\delta\le1$ we have $k\le2A\sqrt K\,\delta^{-s/2}$ and
\[
  1\le mk\delta^s\le 1+k\delta^s.
\]
With $K,A,w$ fixed, $k\delta\to0$ and $k\delta^s\to0$ as $\delta\to0$, since $1<s<2$. We can therefore choose $\delta\le1/2$ sufficiently small to satisfy all the hypotheses of Lemma~\ref{lem:packing} and, in addition,
\begin{equation}\label{eq:normalization}
  1\le mk\delta^s\le1+\varepsilon.
\end{equation}

Enumerate the parent squares by $j=0,\ldots,K-1$. Assign to parent $j$ the $m$ angles
\begin{equation}\label{eq:angles}
  \theta_{j,i}=\theta_0+\frac{w}{Km}(iK+j),
  \qquad i=0,\ldots,m-1.
\end{equation}
Within each parent the angles have spacing $w/m$, so Lemma~\ref{lem:packing} applies with the corresponding initial angle. Across all the parents, the angles form a single angular grid with spacing $w/(Km)$. Every angle in the sector, including its endpoints, lies within $w/(Km)$ of an assigned angle. By \eqref{eq:km},
\[
  \frac{Km}{k}\ge\frac{K\delta^{-s}}{k^2}\ge\frac{1}{4A^2}.
\]
The angular distance is thus at most $4wA^2/k\le4\pi/k$.

Apply Lemma~\ref{lem:packing} inside each parent by translation and scaling by $R$. Each parent now contains exactly $mk$ disjoint child squares of side length $\Delta=R\delta$; assign mass $1/(Kmk)$ to each child. Let $e$ be a direction whose perpendicular, taken modulo $\pi$, lies in the sector. By the angular coverage established above, there is a group of $k$ children lying in the $\sqrt2\Delta$-neighbourhood of a segment of length $RL=4Rk\delta$ making an angle at most $4\pi/k$ with the kernel of $\pi_e$. Their union therefore projects into an interval $J$ of length $D\Delta$, where we may take the absolute constant $D=16\pi+4$: the projected length of the segment is at most $16\pi\Delta$ and the neighbourhood adds at most $2\sqrt2\Delta$. The interval contains the projections of the entire child squares, not just their centres.

Any measure $\mu$ supported on the children and assigning them the prescribed masses satisfies, with $\nu_e=(\pi_e)_\#\mu$,
\[
  \nu_e(J)\ge\frac1{Km}\ge\frac{k\delta^s}{2K}
  \ge\frac{A\delta^{s/2}}{2\sqrt K}.
\]
The second inequality uses $mk\delta^s\le2$. It follows that
\begin{equation}\label{eq:peak}
  \frac{\nu_e(J)}{|J|^{s/2}}
  \ge \frac{A}{2D^{s/2}\sqrt{KR^s}}
  \ge \frac{A}{2D^{s/2}\sqrt{C_0}}.
\end{equation}
This lower bound is preserved under any subsequent refinement that preserves the masses of the children.

Finally, if $K'=Kmk$ and $R'=R\delta$ are the new cardinality and side length, we have
\begin{equation}\label{eq:step}
  K'(R')^s=KR^s(mk\delta^s)\in[KR^s,(1+\varepsilon)KR^s].
\end{equation}
This normalisation allows the step to be repeated without losing the critical exponent.

\section{Iteration and proof of Theorem~\ref{thm:main}}
Start with the square $[0,1]^2$, with $K_0=R_0=1$. Organise the successive refinements into epochs indexed by integers $a\ge1$. In epoch $a$ we divide the circle of segment directions modulo $\pi$ into $a^2$ consecutive sectors of length $\pi/a^2$. Perform one refinement step for each sector, with $A=a$. Number the steps consecutively by $n\ge1$ and use $\varepsilon=2^{-n}$ at step $n$. Denote by $K_n,R_n$ the number and side length of the squares after that step. By \eqref{eq:step},
\begin{equation}\label{eq:invariant}
  1\le K_nR_n^s\le\prod_{j=1}^n(1+2^{-j})\le\exp(1).
\end{equation}
In particular, \eqref{eq:invariant0} holds throughout with $C_0=\exp(1)$. At each step we choose the relative scale $\delta_n\le1/2$ after fixing the parent and sector parameters. We do not require the scales to be the same within an epoch.

Let $E_n$ be the union of the closed squares at level $n$ and let $E=\bigcap_{n\ge0} E_n$. Define a probability measure $\mu$ on $E$ by choosing uniformly among the children of the corresponding square at each step of an infinite path in the tree and pushing forward the path measure to the unique point of intersection of the squares along the path. This point exists and is unique because the squares are nested and $R_n\to0$. The disjointness of the squares ensures that
\begin{equation}\label{eq:cylinder}
 \mu(Q)=1/K_n\qquad\text{for every square }Q\text{ at level }n.
\end{equation}

\subsection*{The Frostman condition}
For $0<r<1$, choose $n\ge1$ such that $R_n\le r<R_{n-1}$. A ball $B(x,r)$ meets only a bounded number of dyadic parent squares of side length $R_{n-1}$, with an absolute bound independent of the generation. Within each parent, the rescaled packing bound \eqref{eq:packing} applies at the relative radius $r/R_{n-1}\in[\delta_n,1)$ and shows that the ball meets at most $C(r/R_n)^s$ children. By \eqref{eq:cylinder} and \eqref{eq:invariant},
\[
  \mu(B(x,r))\lesssim\frac{(r/R_n)^s}{K_n}\le r^s.
\]
For $r\ge1$ it suffices to use $\mu(\mathbb{R}^2)=1$. Thus $\mu$ is $s$-Frostman with a constant depending only on $s$. In particular, the local packing constants do not multiply across generations.

\subsection*{Projected masses}
Fix $e\in S^1$. In each epoch $a$, one of the sectors contains the kernel direction of $\pi_e$. At the corresponding step, \eqref{eq:peak} and \eqref{eq:invariant} yield an interval $J_{a,e}$ such that
\begin{equation}\label{eq:allpeaks}
  (\pi_e)_\#\mu(J_{a,e})\gtrsim a\,|J_{a,e}|^{s/2}.
\end{equation}
The bound holds for the limit measure by \eqref{eq:cylinder}, since the descendants remain inside the selected children. Letting $a\to\infty$ shows that the normalised projected masses for $\mu$ are unbounded. Moreover, the lengths of the intervals tend to zero.

To deduce failure of weak endpoint integrability, set $p_*=2/(2-s)$. If a projected measure has a density $f\in L^{p_*,\infty}(\mathbb{R})$, its decreasing rearrangement $f^*$ satisfies
\[
 f^*(t)\le\|f\|_{L^{p_*,\infty}}t^{-1/p_*},\qquad t>0,
\]
where the usual weak quasi-norm is $\|f\|_{L^{p_*,\infty}}=\sup_{t>0}t^{1/p_*}f^*(t)$. Thus, for every bounded interval $J$,
\begin{equation}\label{eq:weak}
\begin{aligned}
 \int_J |f(x)|\,\mathrm{d}x
 &\le\int_0^{|J|}f^*(t)\,\mathrm{d}t\\
 &\le\frac{p_*}{p_*-1}\|f\|_{L^{p_*,\infty}}|J|^{1-1/p_*}.
\end{aligned}
\end{equation}
Since $1-1/p_*=s/2$, the lower bound \eqref{eq:allpeaks} contradicts \eqref{eq:weak} as $a\to\infty$. This conclusion does not require the projection to be absolutely continuous in every direction: a projection without a density already fails to belong to the indicated space.

\subsection*{Hausdorff measure of the limit set}
It remains to transfer the obstruction to $\mathcal{H}^s|_E$. The Frostman condition and the mass distribution principle give $\mu\lesssim\mathcal{H}^s|_E$ as measures: it suffices to apply the ball estimate to the sets in any cover of small diameter. On the other hand, covering $E$ by the level $N$ squares, we obtain from \eqref{eq:invariant}
\[
  \mathcal{H}^s(E)\le\liminf_{N\to\infty}K_N(\sqrt2 R_N)^s
  \le 2^{s/2}\exp(1),
\]
where we use the convention of defining Hausdorff measure by sums of diameters raised to the power $s$. In particular, $0<\mathcal{H}^s(E)<\infty$.

For a square $Q$ at level $n$, its $K_N/K_n$ descendants at level $N>n$ cover $E\cap Q$. By \eqref{eq:invariant} and \eqref{eq:cylinder},
\begin{align*}
  \mathcal{H}^s(E\cap Q)
  &\le\liminf_{N\to\infty}\frac{K_N}{K_n}(\sqrt2 R_N)^s\\
  &\le\frac{2^{s/2}\exp(1)}{K_n}
  =2^{s/2}\exp(1)\mu(Q).
\end{align*}
The sets $E\cap Q$ form a countable basis of clopen cylinders for the relative topology of $E$: the squares at each level are separated, and their diameters tend to zero. Every relatively open subset of $E$ is a countable disjoint union of such cylinders, obtained by selecting, along each path, the first cylinder contained in the open set. By countable additivity, the preceding bound extends to relatively open sets and, by outer regularity of $\mu$ and monotonicity of $\mathcal{H}^s$, to Borel subsets of $E$. We conclude that
\[
 \mathcal{H}^s|_E\asymp\mu
\]
as measures. In particular, $\lambda=\mathcal{H}^s|_E$ is $s$-Frostman and $(\pi_e)_\#\lambda\gtrsim(\pi_e)_\#\mu$ in every direction. The ratios in \eqref{eq:allpeaks} yield \eqref{eq:obstruction} for $\lambda$. Failure of weak endpoint integrability then follows from \eqref{eq:weak}, and normalisation by $\mathcal{H}^s(E)$ does not affect either conclusion. This completes the proof of Theorem~\ref{thm:main}.\qed

\printbibliography[title={References}]
\end{document}